\documentclass [a4paper,draft,twoside,11pt]{article}

\usepackage[latin1]{inputenc}
\usepackage{amsmath,amsfonts,amssymb}
\usepackage{vmargin,graphicx,theorem}
\usepackage[english]{babel}
\usepackage{enumerate}

\setpapersize[portrait]{A4}
\setmarginsrb{2.2cm}{2cm}{2.2cm}{1cm}{0cm}{0.1cm}{0.5cm}{2cm}

\newcommand{\disp}{\displaystyle}

\newcommand{\dP}{\ensuremath{\mathbb{P}}}

\newcommand{\dR}{\ensuremath{\mathbb{R}}}

\newtheorem{ethm}{Theorem}[section]

\newtheorem{ecor}[ethm]{Corollary}

\newtheorem{eprop}[ethm]{Proposition}

\newtheorem{elem}[ethm]{Lemma}

\newtheorem{edefi}[ethm]{Definition}

\newtheorem{erem}[ethm]{Remark}

\newcommand{\proofend}{~$\rhd$}
\newcommand{\proofbegin}{~$\lhd$}

\newenvironment{eproof}
               {\noindent {\emph{\textbf{Proof}}}\\\proofbegin~}
               {\proofend\\}

\newcommand{\p}[4]{{#3}\!\left#1{#4}\right#2}

\newcommand{\ABS}[1]{\ensuremath{{\left| #1 \right|}}} 
\newcommand{\PAR}[1]{\ensuremath{{\left(#1\right)}}} 
\newcommand{\BRA}[1]{\ensuremath{{\left\{#1\right\}}}} 
\newcommand{\NRM}[1]{\ensuremath{{\left\Vert #1\right\Vert}}} 
\newcommand{\abs}[1]{\ensuremath{\lvert #1 \rvert}}

\newcommand{\nlip}[1]{\ensuremath{\|#1\|_{Lip}}}

\renewcommand{\phi}{\varphi}

\renewcommand{\leq}{\leqslant}
\renewcommand{\geq}{\geqslant}

\newcommand{\varf}[1]{\mathbf{Var}_{#1}}
\newcommand{\entf}[1]{\mathbf{Ent}_{#1}}

\newcommand{\ent}[2]{\p(){\entf{#1}}{#2}}

\newcommand{\var}[2]{\p(){\varf{#1}}{#2}}

\def\disp{\displaystyle}

\newcommand{\R}{\dR}

\newcommand{\C}[1]{\ensuremath{{\mathcal C}^{#1}}}

\newcommand{\Ga}{\boldsymbol{\Gamma}}
\newcommand{\gu}{\boldsymbol{\Ga}}
\newcommand{\gd}{{\bf \gu {\!\!_2}}}

\newcommand{\al}{\alpha}

\newcommand{\be}{\beta}

\newcommand{\ep}{\epsilon}

\renewcommand{\k}{\kappa}
\newcommand{\la}{\lambda}

\newcommand{\e}{\varepsilon}

\begin{document}

\title{\sl
Modified logarithmic Sobolev inequalities in  null curvature}
\author{Ivan Gentil, Arnaud Guillin and Laurent Miclo}

\date{\today}
\maketitle\thispagestyle{empty}

\begin{abstract}
We present a logarithmic Sobolev inequality adapted to a
log-concave measure. Assume that $\Phi$ is a symmetric convex
function on $\dR$ satisfying $(1+\e)\Phi(x)\leq
{x}\Phi'(x)\leq(2-\e)\Phi(x)$ for $x\geq0$ large enough and with
$\e\in]0,1/2]$. We prove that the probability measure on $\dR$
$\mu_\Phi(dx)=e^{-\Phi(x)}/Z_\Phi dx$ satisfies a modified and
adapted logarithmic Sobolev inequality
: there exist three constant $A,B,D>0$ such that
 for all smooth $f>0$,
\begin{equation*}
 \ent{\mu_\Phi}{f^2}\leq
A\int H_{\Phi}\PAR{{\frac{f'}{f}}}f^2d\mu_\Phi, \text{ with }
H_{\Phi}(x)=\left\{
\begin{array}{rl}
\Phi^*\PAR{Bx} &\text{ if }\ABS{x}\geq D,\\
x^2 &\text{ if }\ABS{x}\leq D.
\end{array}
\right.
\end{equation*}
\end{abstract}

\bigskip


\section{Introduction}
\label{sec-in}

A probability measure $\mu$ on
 $\dR^n$ satisfies a logarithmic Sobolev inequality if there exists $C\geq0$ such that,  for every
 smooth enough functions $f$ on $\dR^n$,
\begin{equation}
\label{eq-gross}
\ent{\mu}{f^2} \leq  C \int \ABS{\nabla f}^2 d\mu,
\end{equation}
where
$$ \ent{\mu}{f^2} := \int f^2 \log f^2 d\mu- \int f^2 d\mu  \log \int  f^2 d\mu $$
and where $\ABS{\nabla f}$ is the Euclidean length of the gradient $\nabla f$
of $f$.

\medskip

Gross in~\cite{gross} defines this inequality and shows   that
the canonical Gaussian measure with density $(2\pi )^{-n/2} {e}^{-|x|^2/2}$
with respect to the Lebesgue measure on $\dR^n$ is the basic example of measure
$\mu $ satisfying~\eqref{eq-gross} with the optimal constant  $C=2$.
Since then, many results have presented measures
 satisfying an such  inequality, among them the famous Bakry-\'Emery $\gd$-criterion,
  that we recall now in our particular case.
 Let $\mu(dx)=\exp\PAR{-f(x)}dx$, a probability measure on $\dR^n$ and  assume that there exists $\la>0 $ such that,
\begin{equation}
\label{eq-ba} \forall x\in\dR^n,\quad\text{Hess}\PAR{f(x)}\geq
\la\text{Id},
\end{equation}
in the sense of symmetric matrix. Then Bakry and  \'Emery prove
that $\mu$ is satisfying  inequality~\eqref{eq-gross} with a
optimal constant $0\leq C\leq 2/\la$. We refer to
\cite{bakry-emery,bakry} for the $\gd$-criterion and to
\cite{logsob,ledoux} for a review on logarithmic Sobolev
inequality.

\medskip

The interest of this paper is to give a logarithmic Sobolev
inequality when the probability measure $\mu$ on $\dR$ does'nt
satisfies~\eqref{eq-gross} but it  is still log-concave function
which mean that
 $f''(x)\geq0$.

\medskip

 An answer can be given for
the following measure: Let
 $\alpha\geq 1$ and define the probability measure $\mu_\al$ on
$\dR$ by
\begin{equation}
\label{eq-al}
\mu_\al(dx)=\frac{1}{Z_\al}e^{-\ABS{x}^\al}dx,
\end{equation}
where $Z_\al=\int e^{-\ABS{x}^\al}dx$.

The authors prove, in \cite{ge-gu-mi},  that
 for $1<\alpha<2$,
the measure $\mu_\al$ satisfies the following inequalities,
for all smooth function such that $f\geq0$ and $\int f^2d\mu_\al=1$,
\begin{equation}
\label{logsob1}
\ent{\mu_\al}{f^2}\leq A\var{\mu_\al}{f}+
B\int_{{f}\geq 2}\ABS{\frac{f'}{f}}^\be f^2d\mu_\al,
\end{equation}
where  $A$ and $B$ are some constants  and
$$\var{\mu_\al}{f}:=\int f^2 d\mu_\al-\PAR{\int f d\mu_\al}^2.$$

It is well-known that the probability measure $\mu_\al$ satisfies
 (still for $\alpha\geq 1$)
 a Poincar\'e inequality (or spectral gap inequality) which is  for every smooth enough function $f$,
\begin{equation}
\label{eq-point}
\var{\mu_\al}{f}\leq C\int \ABS{\nabla f}^2d\mu_{\al},
\end{equation}
where  $0<C<\infty$.

Then using~\eqref{eq-point} and~\eqref{logsob1} we get that $\mu_\al$ satisfies
 also this modified logarithmic Sobolev inequality for all smooth
and positive
function  $f$,
\begin{equation}
\label{logsob2}
\ent{\mu}{f^2}\leq C \int H_{a,\al}\PAR{\frac{f'}{f}}f^2 d\mu,
\end{equation}
 here and in the whole paper the convention that $0\cdot \infty=0$ is assumed,
 otherwise stated
where  $a$ and $C$ are positive constants and
$$
H_{a,\al}(x)= \left\{
\begin{array}{ll}
\disp {x^2} &\text{ if } \abs{x}<a, \\
\disp {\abs{x}}^\be\quad &\text{ if } \abs{x}\geq a,
\end{array}
\right.
$$
with  $1/\al+1/\be=1$. The last version of logarithmic Sobolev inequality admits
a $n$ dimensional version, for all smooth function  $f$
on $\dR^n$,
\begin{equation}
\label{logsobn}
\ent{\mu_\al^{\otimes n}}{f^2}\leq C \int H_{a,\al}\PAR{\frac{\nabla f}{f}}f^2 d\mu_\al^{\otimes n},
\end{equation}
where
 by definition we have taken
 \begin{equation}
\label{hs}
  H_{a,\al}\PAR{\frac{\nabla f}{f}}:
 =\sum_{i=1}^n H_{a,\al}\PAR{\frac{\partial_i f}{f}}.
 \end{equation}
Note that Bobkov and Ledoux give in \cite{bobkov-ledoux2}
 a corresponding
result for the critical (exponential) case, when $\al=1$.

\bigskip

Our main purpose here will be to establish the generalization of
inequalities~\eqref{logsob1}, \eqref{logsob2} and \eqref{logsobn}
when the measure on $\dR$ is only a log-concave measure between
$e^{-\ABS{x}}$ and $e^{-x^2}$. More precisely, let $\Phi$ be a
$\C{2}$ convex function on $\dR$. Suppose for simplicity that
$\Phi$ is symmetric. We assume  that $\Phi$ satisfies the
following property, there exists
 $M> 0$
 and $0<\e\leq1/2$ such that $\Phi(M)>0$ and
\begin{equation}
\tag{\bf H} \forall{x\geq M},\,\,\,\,(1+\e)\Phi(x)\leq
x\Phi'(x)\leq(2-\e)\Phi(x)
\end{equation}
We assume during the article that the function $\Phi$ on $\dR$ is satisfying hypothesis $(\bf H)$.

\begin{erem}
\label{rem-int} The assumption $(\bf H)$ implies that there exists
$m_1,m_2>0$ such that
$$
\forall x\geq M,\,\,\,\,m_1x^{1/(1-\e)}\leq\Phi(x)\leq m_2
x^{2-\e}.
$$
This remark explains how, under the hypothesis ${(\bf H)}$, the
function $\Phi$ is between $e^{-\ABS{x}}$ and $e^{-x^2}$.
\end{erem}

Due to the remark~\ref{rem-int}, $\int e^{-\Phi(x)}dx<\infty$.
Then we define the probability measure $\mu_\Phi$ on $\dR$ by
$$
\mu_{\Phi}(dx)=\frac{1}{Z_\Phi}e^{-\Phi(x)}dx,
$$
where $Z_\Phi=\int e^{-\Phi(x)}dx$.

\medskip
The main result of this article is
the following theorem:

\begin{ethm}
\label{theo-ls} Let $\Phi$ satisfying the property $(\bf H)$ then
there exists constants
$A,A',B,D,
\k\geq 0$ such that for any smooth functions
$f\geq0$ satisfying $\int f^2 d\mu_\Phi=1$ we have
\begin{equation}
\label{eq-theo}
 \ent{\mu_\Phi}{f^2}\leq
A\var{\mu_\Phi}{f}+A'\int_{f^2\geq
\k}H_{\Phi}\PAR{{\frac{f'}{f}}}f^2d\mu_\Phi,
\end{equation}
where
\begin{equation}
\label{defh} H_{\Phi}(x)=\left\{
\begin{array}{rl}
\Phi^*\PAR{Bx} &\text{ if }\ABS{x}> D,\\
x^2 &\text{ if }\ABS{x}\leq D,
\end{array}
\right.
\end{equation}
where $\Phi^*$ is the Legendre-Frenchel transform of $\Phi$, $\Phi^*(x):=\sup_{y\in\dR}\BRA{x\cdot y-\Phi(y)}$.
\end{ethm}

It is well known that the measure $\mu_\Phi$  satisfies a Poincar\'e inequality
(inequality~\eqref{eq-point} for the measure $\mu_\Phi$, see for example Chapter~6 of~\cite{logsob}).
 Then we obtain the following corollary:
\begin{ecor}
\label{cor-ls} Let $\Phi$ satisfying the property ${(\bf H)}$ then
there exists $A,B,D\geq0$ such that for any smooth functions $f>0$
 we have
\begin{equation}
\label{eq-cor}
 \ent{\mu_\Phi}{f^2}\leq
A\int H_{\Phi}\PAR{{\frac{f'}{f}}}f^2d\mu_\Phi,
\end{equation}
where $H_\Phi$ is defined on~\eqref{defh}.
\end{ecor}

In~\cite{ge-gu-mi} we investigate some particular example,  where
$\Phi(x)=|x|^\al \log^\beta |x|$, for $\al\in]1,2[$ and
$\be\in\dR$. Theorem~\ref{theo-ls} gives the result in the general
case.

\begin{edefi}
\label{defl} Let $\mu$ a probability measure on $\dR^n$. We said
that $\mu$ satisfies a Logarithmic Sobolev Inequality ($LSI$) of
function $H_\Phi$ $($defined on~\eqref{defh}$)$   if there exists
$A\geq 0$ such that for any smooth functions $f>0$ we have
\begin{equation}
\tag{$LSI$}
 \ent{\mu_\Phi}{f^2}\leq
A\int H_{\Phi}\PAR{{\frac{\nabla f}{f}}}f^2d\mu_\Phi,
\end{equation}
where
 $
 H_{a,\al}\PAR{\frac{\nabla f}{f}}
 $
is defined on~\eqref{hs}.
\end{edefi}
The $LSI$ of function $H_\Phi$ is the $n$-dimensional version of inequality~\eqref{eq-cor}.

\bigskip

In Section~\ref{sec-lsi} we will give the proof of
Theorem~\ref{theo-ls}. It is an adaptation of particular
case studied in~\cite{ge-gu-mi} but it is more technical and
complicated. The proof is cut into two parts, Proposition~\ref{prop-grand2} and~\ref{prop-petit}. In
Subsection~\ref{sec-le}, we will describe the case where the
entropy is large and in
Subsection~\ref{sec-se} we will
study
the other case, when the entropy is small. The two cases are very
different as we can see in the next section but they are connected
to the Hardy's inequality, that we will point out now.

\medskip

Let $\mu,\nu$ be Borel measures on $\R^+$.
 Then the best constant $A$ so that every smooth function
$f$ satisfies
\begin{equation}
\label{eq-hardy}
\int_0^\infty \PAR{f(x)-f(0)}^2d\mu(x)\leq A\int_0^\infty f'^2d\nu
\end{equation}
is finite if and only if
\begin{equation}
\label{eq-cond}
B=\sup_{x>0} \BRA{\mu\PAR{[x,\infty[} \int_0^x\PAR{\frac{d\nu^{ac}}{d t}}^{-1}dt}
\end{equation}
is finite, where $\nu^{ac}$ is the absolutely continuous part of $\nu$ with respect to $\mu$.
Moreover, we have (even if $A$ or $B$ is infinite),
$$B\leq A\leq 4B.$$

One can see for example \cite{bobkov-gotze,logsob} for a review in this domain.

\bigskip

In Section~\ref{sec-app} we will explain some classical properties
of this particular logarithmic Sobolev inequality. We explain
briefly how, as in the classical logarithmic Sobolev inequality of
Gross,
\begin{itemize}
\item The $LSI$ of function $H_\Phi$ satisfies the tensorisation and the perturbation properties.
\item The $LSI$ of function $H_\Phi$ implies also Poincar\'e inequality.
\end{itemize}
The last application proposed is the concentration property for
probability measure satisfying inequality~\eqref{eq-cor}. We
obtain Hoeffding's type inequality, assume that a measure $\mu$ on
$\dR$ satisfies inequality~\eqref{eq-cor} and let $f$ be a Lipschitz
function on $\dR$ with $\nlip{f}\leq1$. Then we get, for some
constants $A,B,D\geq0$ independent of the dimension $n$,
\begin{equation}
\label{eq-con1}
\dP\left({1\over n}\ABS{\sum_{k=1}^nf(X_k)-\mu(f)}>\la\right)\leq
\left\{
\begin{array}{ll}
\disp 2\exp\PAR{-nA \Phi(B\la)}&\text{if }
\la\geq D,\\
\disp 2\exp\PAR{-nA{\la^2}}& \text{if }0\leq\la\leq D,
\end{array}
\right.
\end{equation}
 or equivalently,
\begin{equation}
\label{eq-con2}
\dP\left({1\over \sqrt{n}}\ABS{\sum_{k=1}^nf(X_k)-\mu(f)}>\la\right)\leq
\left\{
\begin{array}{ll}
\disp 2\exp\PAR{-nA \Phi\PAR{B\frac{\la}{\sqrt{n}}}}&\text{if }
\la\geq D\sqrt{n},\\
\disp 2\exp\PAR{-A{\la^2}}& \text{if }0\leq\la\leq D\sqrt{n}.
\end{array}
\right.
\end{equation}
Inequality~\eqref{eq-con2} is interesting because for large enough
$n$ we find the Gaussian concentration, this is natural due to the
convergence of $\frac{1}{\sqrt{n}}\PAR{\sum_{k=1}^nf(X_k)-\mu(f)}$
to the  Gaussian. This result is not a new one, Talagrand
explains  it in~\cite{talagrand95}, see also~\cite{ledoux01} for a
large  review on this topic.

\medskip

Note to finish the introduction that Barthe,
Cattiaux and Roberto~\cite{ca-ba-ro} are studing  the
same sort of log-concave measure. They prove also functional
inequalities with an other point of view, namely Beckner type inequalities or $\Phi$-Sobolev inequalities, in particular one of
their results is concentration inequalities for the same measure
$\mu_\Phi$. Let us also mention that the first author in \cite{G05}, via Prekopa-Leindler inequality, recovers partly our large entropy result.


\section{Proof of logarithmic Sobolev inequality (Theorem~\ref{theo-ls})}
\label{sec-lsi}


Before explaining the proof of Theorem~\ref{theo-ls} we
give a lemma for classical properties satisfied by
the function $\Phi$.

\begin{elem}
\label{lem-av}Assume that $\Phi$ satisfies assumption $(\bf H)$ then
there exists $C\geq 0$ such that for large enough $x\geq 0$,
\begin{equation}
\label{eq-lem1}
x^2\leq C\Phi^*(x),
\end{equation}
\begin{equation}
\label{eq-*}
\e  \Phi\PAR{\Phi'^{-1}(x)}\leq \Phi^*(x)\leq\PAR{1-\e}\Phi\PAR{\Phi'^{-1}(x)},
\end{equation}
\begin{equation}
\label{eq-lem3}
\frac{1}{C}{\Phi'}^{-1}(x)\leq\frac{\Phi^*(x)}{x}\leq C{\Phi'}^{-1}(x).
\end{equation}

\end{elem}
The proof of Lemma~\ref{lem-av} is an easy consequence of the property $(\bf H)$.

\bigskip

For this we will note by {\it smooth function} a
locally
absolutely continuous function on $\dR$. This is the regularity
needed for the use of Hardy inequality in our case.

\subsection{Large entropy}
\label{sec-le}


The proof of $LSI$ for large entropy is based on the next lemma,
we give a $LSI$ {\it saturate on the left}.

\begin{elem}
\label{lem-h}
Let $h$ defined as follows

\begin{equation}
\label{def-h}
h(x)=\left\{
\begin{array}{lr}
1&\text{if}\,\,\ABS{x}< M\\
\frac{x^{2}}{\Phi(x)}&\text{if}\,\,\ABS{x}\geq M.
\end{array}
\right.
\end{equation}
Then there exists ${C_h}\geq0$ such that  for every smooth
function $g$ we have
\begin{equation}
\label{eq-lemh}
\ent{\mu_\Phi}{g^2}\leq {C_h}\int g'^2 hd\mu_\Phi.
\end{equation}
\end{elem}
\begin{eproof}
We use Theorem~3 of \cite{barthe-roberto} which is a refinement of the criterion of a Bobkov-G\"otze theorem
(see Theorem~5.3 of \cite{bobkov-gotze}).

The constant ${C_h}$ satisfies $\max(b_{-},b_{+})\leq C_h\leq \max(B_{-},B_{+})$ where
$$
b_{+}=\sup_{x\geq0}
\mu_\Phi([x,+\infty[)\log\PAR{1+\frac{1}{2\mu_\Phi([x,+\infty[)}}
\int_0^x  {Z_\Phi}\frac{e^{\Phi(t)}}{h(t)}dt,
$$
$$
b_{-}=\sup_{x\leq0}
\mu_\Phi(]-\infty,x])\log\PAR{1+\frac{1}{2\mu_\Phi(]-\infty,x])}}
\int_x^0  {Z_\Phi}\frac{e^{\Phi(t)}}{h(t)}dt,
$$
$$
B_{+}=\sup_{x\geq0}
\mu_\Phi([x,+\infty[)\log\PAR{1+\frac{e^2}{\mu_\Phi([x,+\infty[)}}
\int_0^x  {Z_\Phi}\frac{e^{\Phi(t)}}{h(t)}dt,
$$
$$
B_{-}=\sup_{x\leq0}
\mu_\Phi(]-\infty,x])\log\PAR{1+\frac{e^2}{\mu_\Phi([-\infty,x[)}}
\int_x^0  {Z_\Phi}\frac{e^{\Phi(t)}}{h(t)}dt.
$$

An easy approximation proves that for large positive $x$
\begin{equation}
\label{eq-approx}
\mu_\Phi([x,\infty[)=\int_x^\infty \frac{1}{Z_\Phi}e^{-\Phi(t)}dt
\sim_\infty \frac{1}{Z_\Phi \Phi'(x)}{e^{-\Phi(x)}},
\end{equation}
and
$$
\int_0^x  {Z_\Phi}\frac{e^{\Phi(t)}}{h(t)}dt\sim_\infty \frac{Z_\Phi}{h(x)\Phi'(x)}e^{\Phi(x)},
$$
and one may prove similar behaviors for negative $x$.

Then, there is $K$ such that for $x\geq M$,
$$
\mu_\Phi([x,+\infty[)\log\PAR{1+\frac{1}{2\mu_\Phi([x,+\infty[)}}
\int_0^x  {Z_\Phi}\frac{e^{\Phi(t)}}{h(t)}dt\leq K\frac{\Phi(x)}{\Phi'(x)^2h(x)}=
K\PAR{\frac{\Phi(x)}{x\Phi'(x)}}^2.
$$
The right hand term is bounded by the assumption $({\bf H})$.

A simple calculation then yields that
constants $b_{+}$, $b_{-}$, $B_{+}$ and  $B_{-}$ are finite and the lemma is proved.
\end{eproof}

\begin{erem}
Note that this lemma can be proved in a more general case, when
$\Phi$ does not satisfy hypothesis $({\bf H})$. In
\cite{bobkov-ledoux1} the authors prove this result for the
symmetric exponential measure.
\end{erem}


\begin{eprop}
\label{prop-grand2} There exists $A,B,D,A'\geq0$ such that for
any functions $f\geq0$ satisfying
$$\int f^2 d\mu_\Phi=1\,\,\text{and}\,\,\ent{\mu_\Phi}{f^2}\geq1$$
 we have
\begin{equation}
\label{eq-grand} \ent{\mu_\Phi}{f^2}\leq
A'\var{\mu_\Phi}{f}+A\int_{f\geq2}H_\Phi\PAR{{\frac{f'}{f}}}d\mu_\Phi,
\end{equation}
where
\begin{equation*}
H_\Phi(x)=\left\{
\begin{array}{rl}
\Phi^*\PAR{Bx} &\text{ if }\ABS{x}\geq D,\\
x^2 &\text{ if }\ABS{x}\leq D.
\end{array}
\right.
\end{equation*}
\end{eprop}

As we will see in the proof, $A'$ does not depend on the function
$\Phi$.

{\noindent {\emph{\textbf{Proof of Proposition~\ref{prop-grand2}}}}\\~}
\proofbegin~
Let $f\geq0$ satisfying $\int f^2d\mu_\Phi=1$.

A careful study of the function
$$
x\rightarrow -x^2\log x^2 +5(x-1)^2+x^2-1+(x-2)_+^2\log (x-2)_+^2
$$
proves that for every $x\geq 0$
$$
x^2\log x^2\leq 5(x-1)^2+x^2-1+(x-2)_+^2\log (x-2)_+^2.
$$

We know that $\int (f-1)^2d\mu_\Phi\leq 2\var{\mu_\Phi}{f}$,
recalling that $\int f^2d\mu_\Phi=1$ and $f\geq0$,
\begin{eqnarray*}
\int f^2\log f^2d\mu_\Phi & \leq &5\int (f-1)^2d\mu_\Phi +\int (f^2-1)d\mu_\Phi+
\int (f-2)_+^2\log (f-2)_+^2d\mu_\Phi\\
 & \leq &10\var{\mu_\Phi}{f}+ \int (f-2)_+^2\log (f-2)_+^2d\mu_\Phi.
\end{eqnarray*}

Since $\int f^2 d\mu_\Phi=1$, one can easily prove that
$$
\int (f-2)_+^2d\mu_\Phi\leq 1,
$$
then $\int (f-2)_+^2\log (f-2)_+^2d\mu_\Phi\leq \ent{\mu_\Phi}{(f-2)_+^2},$
and
\begin{eqnarray*}
\ent{\mu_\Phi}{f^2}\leq10\var{\mu_\Phi}{f}+\ent{\mu_\Phi}{(f-2)_+^2}.
\end{eqnarray*}

Hardy's inequality of Lemma~\ref{lem-h} with $g=(f-2)_+$ gives

\begin{equation}
\label{eq-br}
\ent{\mu_\Phi}{(f-2)_+^2}\leq {C_h}\int (f-2)_+'^2 hd\mu_\Phi
={C_h}\int_{f\geq2} f'^2 hd\mu_\Phi.
\end{equation}

Due to the assumption $(\bf H)$, the function
$
h(x)=x^2/\Phi(x),
$
is increasing on $[M,\infty[$ and
$$
\lim_{x\rightarrow\infty}h(x)=\infty.
$$
We can assume that $\Phi(M)>0$. We note $m=h(M)>0$
Let us define the function $\tau$ as follow
\begin{equation}
\label{def-tau}
\tau(x)=\left\{
\begin{array}{lr}
x\Phi\PAR{h^{-1}(m)}/(8C_hm)&\text{if}\,\,0\leq x\leq m\\
\Phi\PAR{h^{-1}(x)}/(8C_h)&\text{if}\,\,x\geq m
\end{array}
\right.
\end{equation}

For all $x\geq M$, we have
$\tau(h(x))=\Phi(x)/(8C_h)$ and then, an easy calculus gives
that $\tau$ is increasing on $[0,\infty[$.

Let $u>0$,
\begin{multline*}
C_h\int_{f\geq2} f'^2 hd\mu_\Phi=C_h\int_{f\geq2}
u\PAR{\frac{f'}{f}}^2 \frac{h}{u}f^2d\mu_\Phi\\
\leq C_h\int_{f\geq2} \tau^*\BRA{u\PAR{\frac{f'}{f}}^2}
f^2d\mu_\Phi+ \int_{f\geq2}C_h\tau\PAR{\frac{h}{u}}f^2d\mu_\Phi
\end{multline*}

For every function $f$ such that $\int f^2d\mu_\Phi=1$ and
for every measurable function $g$ such that $\int f^2gd\mu_\Phi$
exists we get
$$
\int f^2g d\mu\leq\ent{\mu_\Phi}{f^2}+\log{\int e^gd\mu_\Phi}.
$$
Indeed, this inequality is also true for all function $g\geq0$ even
 if the above integrals are infinite.
This inequality is also true for all function $g\geq0$ even
integrals are infinite.

We apply the previous inequality with $g=4C_h\tau\PAR{{h}/{u}}$
and we obtain
$$
\int_{f\geq2}C_h\tau\PAR{\frac{h}{u}}f^2d\mu_\Phi
\leq\frac{1}{4}\int 4C_h\tau\PAR{\frac{h}{u}}f^2d\mu_\Phi
\leq \frac{1}{4}\PAR{\ent{\mu_\Phi}{f^2}+\log\int e^{4C_h\tau\PAR{\frac{h}{u}}}d\mu_\Phi}.
$$
If $u=1$ we have, by construction,  $\int e^{4C_h\tau\PAR{\frac{h}{u}}}d\mu_\Phi<\infty$, then
we get
$$
\lim_{u\rightarrow\infty}\int e^{4C_h\tau\PAR{\frac{h}{u}}}d\mu_\Phi=1.
$$
Then, by the bounded convergence theorem, there exists
$u_0$ such that $\int e^{4C_h\tau\PAR{\frac{h}{u_0}}}d\mu_\Phi\leq e$.

Thus we have
$$\ent{\mu_\Phi}{{(f-2)_+}^2}\leq  C_h\int_{f\geq2} \tau^*\BRA{u_0\PAR{\frac{f'}{f}}^2} f^2d\mu_\Phi
+\frac{1}{4}\ent{\mu_\Phi}{f^2}+\frac{1}{4}$$
$\ent{\mu_\Phi}{f^2}\geq 1$, implies
\begin{equation*}
\ent{\mu_\Phi}{f^2}\leq 20\var{\mu_\Phi}{f}+
2C_h\int_{f\geq2} \tau^*\BRA{u_0\PAR{\frac{f'}{f}}^2} f^2d\mu_\Phi.
\end{equation*}

Then Lemma~\ref{lem-legen} gives the proof of inequality~\eqref{eq-grand}.
\proofend


\begin{elem}
\label{lem-legen}

There exist constants $A,B,C,D\geq0$ such that

\begin{equation*}
\forall x\geq 0,\,\,\,\tau^*(x^2)\leq
\left\{
\begin{array}{lr}
A\Phi^*(Cx)&\text{ if } x\geq D, \\
Bx^2&\text{ if } x\leq D.
\end{array}
\right.
\end{equation*}
\end{elem}

\begin{eproof}
Let $x>0$,
$$
\tau^*(x)=\sup_{y\geq0}\BRA{xy-\tau(y)}.
$$
Let $m=h(M)>0$, then
\begin{equation*}
\begin{array}{rl}
\disp\tau^*(x)&=\disp\max\BRA{\sup_{y\in[0,m[}\BRA{xy-\tau(y)},\sup_{y\geq m}\BRA{xy-\tau(y)}},\\
&\disp\leq \sup_{y\in[0,m[}\BRA{xy-\tau(y)}+\sup_{y\geq m}\BRA{xy-\tau(y)}.
\end{array}
\end{equation*}
We have $\disp\sup_{y\in[0,m[}\BRA{xy-\tau(y)}\leq xm$,
because $\tau$ is positive. Then the definition of $\tau$
 implies that
$$
\sup_{y\geq m}\BRA{xy-\tau(y)}=\sup_{y\geq M}\BRA{x\frac{y^2}{\Phi(y)}-\frac{\Phi(y)}{8C_h}}.$$

Let define $\psi_x(y)=x{y^2}/{\Phi(y)}-{\Phi(y)}/\PAR{8C_h}$ for $y\geq M$. We have
$$
\psi_x'(y)=xy\frac{2\Phi(y)-y\Phi'(y)}{\Phi^2(y)}-\frac{\Phi'(y)}{8C_h}.
$$
Due to the property $(\bf H)$, there is $D>0$ such that
$$
\forall x\geq D,\,\,\,\,\sup_{y\geq m}\BRA{xy-\tau(y)}=x\frac{y^2_x}{\Phi(y_x)}-\frac{\Phi(y_x)}{8C_h},
$$
where $y_x\geq M$ satisfies
$$
x=\frac{1}{8C_h}\frac{\Phi'(y_x)\Phi^2(y_x)}{y_x\PAR{2\Phi(y_x)-y_x\Phi'(y_x)}}.
$$
The assumption $(\bf H)$ implies that
$$
\e y_x\Phi'(y_x)\leq 2\Phi(y_x)-y_x \Phi'(y_x)\leq\frac{1-\e}{1+\e}y_x\Phi'(y_x),
$$
then
\begin{equation}
\label{eq-x}
\frac{1}{8C_h(1-\e)(2-\e)}\Phi'^2(y_x)\leq x\leq\frac{1}{8C_h\e(1+\e)}\Phi'^2(y_x).
\end{equation}
We get with the assumption  $(\bf H)$,
\begin{equation*}
\begin{array}{rl}
\displaystyle \forall x\geq D,\,\,\,\,\sup_{y\geq m}\BRA{xy-\tau(y)}
& \displaystyle\leq \frac{1}{8C_h\e(1+\e)}\Phi'^2(y_x)\frac{y_x^2}{\Phi(y_x)}-\frac{1}{8C_y}\Phi(y_x)\\
&\leq\displaystyle \frac{(2-\e)^2}{8C_h\e(1+\e)}\Phi(y_x).
\end{array}
\end{equation*}
Equation~\eqref{eq-x} gives,
\begin{equation*}
y_x\leq \Phi'^{-1}(\sqrt{Cx})
\end{equation*}
where $C>0$. Then  we get
\begin{equation*}
\begin{array}{rl}
\displaystyle\forall x\geq D,\,\,\,\, \sup_{y\geq m}\BRA{xy-\tau(y)}
& \displaystyle\leq  \frac{(2-\e)^2}{8C_h\e(1+\e)}\Phi\PAR{ \Phi'^{-1}(\sqrt{Cx})}.
\end{array}
\end{equation*}

We obtain, using inequality~\eqref{eq-*} of Lemma~\ref{lem-av},
$$
\forall x\geq D,\,\,\,\,\sup_{y\geq m}\BRA{xy-\tau(y)}\leq\frac{1}{8C_h\e^2(1+\e)}\Phi^*\PAR{\sqrt{Cx}}.
$$
then,
$$
\forall x\geq D,\,\,\,\,\tau^*(x)\leq xm+K\Phi^*\PAR{\sqrt{Cx}}.
$$
Using inequality~\eqref{eq-lem1} of Lemma~\ref{lem-av} we get
$$
\forall x\geq D,\,\,\,
\tau^*(x)\leq K'\Phi^*\PAR{\sqrt{Cx}},
$$
for some $K'\geq0$.

On the other hand, the function $\tau$ is non-negative and satisfy $\tau(0)=0$ then $\tau^*(0)=0$.
$\tau^*$ is also a convex function, then there exists $m'$ such that
$$
\forall x\in[0,D],\,\,\,\,
\tau^*(x)\leq xm',
$$
which proves the lemma.
\end{eproof}

\begin{ecor}
For any smooth function $f>0$ on $\dR$ satisfying
$$
\int f^2d\mu_\Phi=1,\text{ and }\ent{\mu_\Phi}{f^2}\geq1,
$$
we have
\begin{equation*}
\ent{\mu_\Phi}{f^2}\leq C\int{H_{\Phi}\PAR{{\frac{f'}{f}}}
f^2}d\mu_\Phi,
\end{equation*}
where
\begin{equation*}
H_{\Phi}(x)=\left\{
\begin{array}{rl}
\Phi^*\PAR{Bx} &\text{ if }\ABS{x}\geq D\\
x^2 &\text{ if }\ABS{x}\leq D,
\end{array}
\right.
\end{equation*}
and $B,D\geq0$.
\end{ecor}

\begin{eproof}
Due to the property $(\bf H)$ the measure $\mu_\Phi$ satisfies a Spectral Gap inequality,
$$
\var{\mu_\Phi}{f}\leq C_{SG}\int f'^2d\mu_\Phi,
$$
with $C_{SG}\geq0$. We apply inequality~\eqref{eq-grand} to get
the result.
\end{eproof}


\subsection{Small entropy} \label{sec-se}


\begin{elem}
\label{lem-con}
Let $\la>0$ and define  the function $\psi$ by
\begin{equation*}
\psi(x) =\BRA{(\Phi^*)^{-1}\PAR{\la\log x}}^2.
\end{equation*}

Then for all $\la>0$ there exists $A_\la>0$ such that the function
$\psi$ is well defined, positive, increasing, concave on  $[A_\la,\infty[$ and satisfies
$\psi(A_\la)\geq1$.
\end{elem}

\begin{eproof}
Let $\la>0$ be fixed.
 Classical property of the Legendre-Frenchel transform
implies that $\Phi^*$ is convex. Due to the property
($\bf H$),
$\PAR{\Phi^*}^{-1}(\la \log x)$ is well defined for $x\geq M_1$
with $M_1>0$. Then we get  on $[M_1,\infty[$,
$$
\psi'(x)=2g'(\la \log x)g(\la \log x)\frac{\la}{x},
$$
and
$$
\psi''(x)=2 g(\la \log x)\frac{\la^2}{x^2}\PAR{g''(\la\log x)-\frac{g'(\la\log x)}{\la}+\frac{{g'}^2(\la\log x)}{g(\la\log x)}},
$$
where, for simplicity, we have noted  $g=\PAR{\Phi^*}^{-1}$.

For $x$ large enough $g$ is non-negative and increasing  and then $\psi$ is increasing on  $[M_2,\infty[$, with $M_2\geq0$.

An easy estimation gives that
 as $x$ goes to infinity,
\begin{equation}
\label{eq-est}
\frac{{g'}(x)}{g(x)}=o_\infty\PAR{1},
\end{equation}
then since $\PAR{\Phi^*}^{-1}$ is concave, for all large enough
$x$, $\psi''(x)\leq 0$. Then one can find $A_\la>0$ such that
properties on the Lemma~\ref{lem-con} are
true.
\end{eproof}

The proof of $LSI$ for small  entropy is based on the next lemma,
we give a $LSI$ {\it saturate on the right}.

\begin{elem}
\label{lem-g} There exists $\la>0$ which
depends on the function
$\Phi$ such that if we note by $A_\la$ the constant of
Lemma~\ref{lem-con} we get for all $g$ defined on $[T,\infty[$
with $T\in[T_1,T_2]$ for some fixed $T_1,T_2$,
 and verifying that
$$g(T)=\sqrt{A_\la},\,\, g\geq \sqrt{A_\la}\,\, \text{and}\,\, \int_{T}^\infty g^2d\mu_\al\leq 2A_\la+2.$$
Then we get
\begin{equation}
\label{eq-lemme}
\int_{T}^\infty (g-\sqrt{A_\la})_+^2\psi(g^2)\mu_\Phi\leq
 {C_1}\int_{[T,\infty[} {g'}^2d\mu_\Phi,
\end{equation}
where
$\psi$ is defined on Lemma~\ref{lem-con}.

\medskip

The constant ${C_1}$ depend on $\Phi$ and $\la$ but does not depend
on the value of $T\in[T_1,T_2]$.
\end{elem}

\begin{eproof}
Let use Hardy's inequality as explained in the introduction.
We have $g(T)=A_\la$. We apply inequality~\eqref{eq-hardy} on $[T,\infty[$ with the function
$(g-\sqrt{A_\la})_+$
and the following measures
$$d\mu=\psi(g^2)d\mu_\Phi\,\,\text{and}\,\,\nu=\mu_\Phi.$$
Then  the constant $C$ in inequality~\eqref{eq-lemme} is finite if and only if
$$
B=\sup_{x\geq T} \int_{T}^x e^{\Phi(t)}dt\int_x^\infty
\psi(g^2)d\mu_\Phi,
$$
is finite.

By Lemma~\ref{lem-con}, $\psi$ is concave on $[A_\la,\infty[$ then by Jensen inequality, for all $x\geq T$ we get
$$
 \int_x^\infty \psi\PAR{ g^2}d\mu_\Phi
\leq  \mu_\Phi\PAR{[x,\infty[}\psi\PAR{\frac{\int_{x}^\infty {g^2}d\mu_\Phi}
{\mu_{\Phi}\PAR{[x,\infty[}}}.
$$

Then we have
\begin{eqnarray}
\label{eb}
B  \leq  \sup_{x>T_1}\BRA{\int_{T_1}^x  e^{\Phi(t)}dt
\mu_\Phi\PAR{[x,\infty[}\psi\PAR{\frac{\int_{x}^\infty g^2d\mu_\Phi}
{\mu_{\Phi}\PAR{[x,\infty[}}}}
\end{eqnarray}

Due to the property
${(\bf H)}$ there exists $K>1$ such that
\begin{equation}
\label{e-k}
\Phi'(x)e^{\Phi(x)}\leq e^{K\Phi(x)},
\end{equation}
and
$$
\int_{T_{_1}}^xe^{\Phi(t)}dt\leq \frac{Ke^{\Phi(x)}}{\Phi'(x)},\quad
\int_{x}^\infty e^{-\Phi(t)}dt\leq \frac{e^{-\Phi(x)}}{\Phi'(x)},
$$
for large enough $x$. By~\eqref{e-k} we get also for large enough $x$ that
$$
e^{-K\Phi(x)}\leq\int_{x}^\infty e^{-\Phi(t)}dt.
$$

Then for large enough $x$,
 uniformly in the previous $g$, one have
$$
\int_{T_{_1}}^xe^{\Phi(t)}dt\mu_\Phi\PAR{[x,\infty[}
\psi\PAR{\frac{\int_x^\infty g^2d\mu_\Phi}{\mu_\Phi\PAR{[x,\infty[}}}
\leq\frac{K}{\PAR{\Phi'(x)}^2}
\psi\PAR{\frac{\int_x^\infty g^2d\mu_\Phi}{K}e^{K\Phi(x)}}.
$$

For $x$ large enough, $$\frac{\int_x^\infty g^2d\mu_\Phi}{K}\leq1.$$

Then, by definition of  $\psi$, for large enough $x$,
$$
\int_{T_{_1}}^xe^{\Phi(t)}dt\mu_\Phi\PAR{[x,\infty[}
\psi\PAR{\frac{ \int_x^\infty g^2d\mu_\Phi }{\mu_\Phi\PAR{[x,\infty[}}}
\leq K\PAR{\frac{{\Phi^{*}}^{-1}\PAR{\la K \Phi(x)}}{{\Phi'(x)}}}^2.
$$
There is also $C_\ep$ such that, for $x$ large enough
$$
{\Phi^{*}}^{-1}(x)\leq\Phi'\PAR{\Phi^{-1}(C_\ep x)},
$$
as one can see from equation~\eqref{eq-*}.

Then one can choose $\la=1/(KC_\ep)$ and the
lemma is proved. Note that $\la$ depends only on the function $\Phi$.

The constant $B$ on~\eqref{eb} is bounded by $K$
which does'nt depend on $T$ on $[T_1,T_2]$.
\end{eproof}

\begin{eprop}
\label{prop-petit}
There exists $A,A',B,D>0$ such that for any functions $f\geq0$ satisfying
$$\int f^2 d\mu_\al=1\,\, \text{and} \,\,\ent{\mu_\al}{f^2}\leq 1$$
 we have
$$
\ent{\mu_\Phi}{f^2}\leq A\var{\mu_\Phi}{f}+A'\int_{f^2\geq
A_\la}H\PAR{{\frac{f'}{f}}}d\mu_\Phi,
$$

where
\begin{equation*}
H(x)=\left\{
\begin{array}{rl}
\Phi^*\PAR{Bx} &\text{ if }\ABS{x}\geq D,\\
x^2 &\text{ if }\ABS{x}\leq D.
\end{array}
\right.
\end{equation*}
\end{eprop}

\begin{eproof}
Let $f\geq0$ satisfying $\int f^2d\mu_\al=1$.

We can assume that $A_\la\geq 2$.
A careful study of the function
$$
x\rightarrow -x^2\log x^2 +A(x-1)^2+x^2-1+(x-\sqrt{A_\la})_+^2\log (x-\sqrt{A_\la})_+^2
$$
proves that there exists $A$ such that for every $x\in \R^+$
$$
x^2\log x^2\leq A(x-1)^2+x^2-1+(x-\sqrt{A_\la})_+^2\log (x-\sqrt{A_\la})_+^2.
$$

Then we get
\begin{eqnarray}
\label{eq-dep}
\ent{\mu_\al}{f^2}=\int f^2\log f^2d\mu_\al \leq A\var{\mu_\al}{f}+ \int \PAR{f-\sqrt{A_\la}}_+^2\log f^2d\mu_\al,
\end{eqnarray}
where $\sqrt{A_\la}$ is defined
as in
Lemmas~\ref{lem-con} and~\ref{lem-g}.

\bigskip

Fix
$\la$ as in Lemma~\ref{lem-g}. We define the function $K$ on $[A_\la,\infty[$ by
$$
K(x)=\sqrt{\frac{{\log x^2}}{\psi(x^2)}},
$$
where $A_\la$ is defined on Lemma~\ref{lem-g}.

Let now define $T_1<T_2$ such that
$$
\mu_\Phi\PAR{]\infty,T_1]}=\frac{3}{8},\,\,\mu_\Phi\PAR{[T_1,T_2]}=\frac{1}{4}
\,\,\text{and}\,\,\mu_\Phi\PAR{[T_2,+\infty[}=\frac{3}{8}.
$$
Since $\int f^2d\mu_\Phi=1$ there exists $T\in [T_1,T_2]$ such that $f(T)\leq A_\lambda$.

Let us define $g$ on $[T_1,\infty]$ as follow
\begin{equation*}
g= \sqrt{A_\la}+\PAR{f-\sqrt{A_\la}}_+K(f)\text{ on } [T,\infty[.
\end{equation*}
 Function $g$ satisfies $g(T)={\sqrt{A_\la}}$ and $g(x)\geq {\sqrt{A_\la}}$ for all $x\geq T$.

Then we have
\begin{eqnarray}
\begin{array}{rl}
\disp\int_ T^\infty g^2d\mu_\Phi &\disp\leq\int_{T_1}^\infty g^2d\mu_\Phi\\
\disp &\disp\leq  2A_\la+2\int_{[T_1,\infty[\cap\BRA{f^2\geq{A_\la}}} f^2 K^2(f)d\mu_\Phi\\
\disp &\disp\leq  2A_\la+2\int_{[T_1,\infty[} f^2 \log(f^2)d\mu_\Phi\\
\disp &\disp\leq 2A_\la+2,
\end{array}
\end{eqnarray}
where we are using the growth of $\psi$ on $[A_\la,\infty[$ and $\psi(A_\la)\geq 1$.

Assumptions on Lemma~\ref{lem-g} are satisfied, we obtain by inequality~\eqref{eq-lemme}

\begin{equation*}
\int_{T}^\infty (g-\sqrt{A_\lambda})_+^2\psi(g^2)d\mu_\Phi
\leq {C_1}\int_{[T,\infty[\cap\BRA{f^2\geq A_\la}} {g'}^2d\mu_\Phi.
\end{equation*}
Let us compare the various terms now.

Due to the property (${\bf H}$), $K$ is lower bounded on $[\sqrt{A_\la},\infty[$ by
$\al\geq1$ (maybe for $A_\la$ larger), then
we get firstly
$$
\sqrt{
A_\la}+\PAR{f-\sqrt{A_\la}}_+K(f)\geq\sqrt{A_\la}+\PAR{f-\sqrt{A_\la}}_+
\al\geq
f\quad \text{on} \quad\BRA{f^2\geq{A_\la}}.
$$

Then
\begin{multline*}
\PAR{g-\sqrt{A_\la}}_+^2\psi(g^2)=\PAR{f-\sqrt{A_\la}}_+^2
K(f)^2\psi\PAR{\sqrt{A_\lambda}+\PAR{f-\sqrt{A_\la}}_+K(f)}^2 \\
\geq\PAR{f-\sqrt{A_\la}}_+^2K(f)^2\psi(f^2)
=\PAR{f-\sqrt{A_\la}}_+^2\log f^2,
\end{multline*}
by the definition of $K$, then we obtain
\begin{equation}
\label{eq-first}
 \int_{T}^\infty (f-\sqrt{A_\la})_+^2\log f^2d\mu_\Phi
\leq \int_{T}^\infty (g-\sqrt{A_\la})_+^2\psi(g^2)d\mu_\Phi.
\end{equation}

Secondly we have on  $\BRA{f\geq\sqrt{A_\la}}$
\begin{eqnarray*}
\begin{array}{rl}
g' & =\disp f'K(f)+\PAR{f-\sqrt{A_\la}}_+f'K'(f)\\
 & = \disp f'K(f)\PAR{1+ \PAR{f-\sqrt{A_\la}}_+\frac{K'(f)}{K(f)}}
\end{array}
\end{eqnarray*}

 But
we have for $x\geq\sqrt{A_\la}$
\begin{eqnarray*}
\begin{array}{rl}
\disp\ABS{1+(x-\sqrt{A_\la})\frac{K'(x)}{K(x)}}& \disp\leq 1+x\ABS{\frac{K'(x)}{K(x)}}\\
& \leq \disp 1+\frac{1}{2\log x}+
\ABS{\frac{\la}{x}\frac{g'(\la2\log x)}{g(\la2\log x)}},
\end{array}
\end{eqnarray*}
where $g(x)={\Phi^*}^{-1}(x)$. Using Lemma~\ref{lem-con} and the estimation~\eqref{eq-est} we obtain that
there exists $C>0$ such that for all $x\geq \sqrt{A_\la}$,
$$
\ABS{1+(x-\sqrt{A_\la})\frac{K'(x)}{K(x)}}\leq C
$$

We get then
$$
g'^2\leq C f'^2K^2(f)\quad\text{on}\quad\BRA{f^2\geq{A_\la}},
$$
for some $C<\infty$ and then
\begin{equation}
\label{eq-second}
\int_{[T,\infty[ \cap \BRA{f^2\geq {A_\la}}}g'^2d\mu_\Phi
\leq C\int_{[T,\infty[ \cap \BRA{f^2\geq {A_\la}}}f'^2K^2(f)d\mu_\Phi.
\end{equation}

By equation~\eqref{eq-first} and~\eqref{eq-second} we obtain
$$
 \int_{T}^\infty (f-\sqrt{A_\la})_+^2\log f^2d\mu_\Phi\leq C\int_{[T,\infty[ \cap \BRA{f^2\geq {A_\la}}}f'^2K^2(f)d\mu_\Phi.
$$

Let $u_0>0$,
\begin{multline*}
 \int_{T}^\infty (f-\sqrt{A_\la})_+^2\log f^2d\mu_\Phi
\leq \\
C\int_{[T,\infty[ \cap \BRA{f^2\geq {A_\la}}} \tau^*_2\PAR{u_0\PAR{\frac{f'}{f}}^2}f^2d\mu_\Phi+
\int_{[T,\infty[ \cap \BRA{f^2\geq {A_\la}}}\tau_2\PAR{\frac{K^2(f)}{u_0}}f^2d\mu_\Phi,
\end{multline*}
where the function $\tau_2$ is defined as in equation~\eqref{def-tau} by
\begin{equation}
\label{def-tau_2}
\tau_2(x)=\left\{
\begin{array}{lr}
\disp x\Phi\PAR{h^{-1}(m)}\frac{1-\ep}{2\la m}&\text{if}\,\,0\leq x < m\\
\disp \Phi\PAR{h^{-1}(x)}\frac{1-\ep}{2\la}&\text{if}\,\,x\geq m,
\end{array}
\right.
\end{equation}
where $h$ is defined on equation~\eqref{def-h} and $m$ on equation~\eqref{def-tau}. The function $\tau_2$ is equal to $\tau$
up to  a constant factor.

Using Lemma~\ref{lem-tau2}
we get
\begin{multline*}
 \int_{T}^\infty (f-\sqrt{A_\la})_+^2\log f^2d\mu_\Phi
\leq \\
C\int_{[T,\infty[ \cap \BRA{f^2\geq {A_\la}}} \tau_2^*\PAR{u_0\PAR{\frac{f'}{f}}^2}f^2d\mu_\Phi+
\frac{1}{2}\int_{[T,\infty[ \cap \BRA{f^2\geq {A_\la}}}f^2\log f^2d\mu_\Phi.
\end{multline*}

The same method can be used on $]-\infty,T]$ and then there is $C'<\infty$ such that
\begin{multline*}
 \int_{-\infty}^T (f-\sqrt{A_\la})_+^2\log f^2d\mu_\Phi
\leq \\
C'\int_{]-\infty,T] \cap \BRA{f^2\geq {A_\la}}} \tau_2^*\PAR{u_0\PAR{\frac{f'}{f}}^2}f^2d\mu_\Phi+
\frac{1}{2}\int_{]-\infty ,T]\cap \BRA{f^2\geq {A_\la}}}f^2\log f^2d\mu_\Phi.
\end{multline*}
And then we get
\begin{multline*}
 \int (f-\sqrt{A_\la})_+^2\log f^2d\mu_\Phi
\leq
(C+C')\int_{\BRA{f^2\geq {A_\la}}} \tau_2^*\PAR{u_0\PAR{\frac{f'}{f}}^2}f^2d\mu_\Phi+
\frac{1}{2}\int_{\BRA{f^2\geq {A_\la}}}f^2\log f^2d\mu_\Phi.
\end{multline*}
Note that constants $C$ and $C'$ don't depend on $T\in[T_1,T_2]$.

Then by inequality~\eqref{eq-dep} and Lemma~\ref{lem-legen}, Proposition~\ref{prop-petit} is proved.
\end{eproof}

\begin{elem}
\label{lem-tau2}
There exists $u_0>0$ such that, for all $x\geq A_\la$  we have
$$
\tau_2\PAR{\frac{K^2(x)}{u_0}}\leq\frac{1}{2}\log x^2.
$$
\end{elem}

\begin{eproof}
Let $\k =2\la/(1-\ep)$.

For all $x\geq M$, where $M$ is defined on equation~\eqref{def-tau},  we have
$$
\tau_2\PAR{h(x)}=\frac{\Phi(x)}{\k},
$$
$$
\tau_2\PAR{\frac{x^2}{\Phi(x)}}=\frac{\Phi(x)}{\k}.
$$
$\tau_2$ is increasing, then due to the property ($\bf H$) we have for $x\geq M$
$$
\tau_2\PAR{(1+\ep)^2\frac{\Phi(x)}{\Phi'(x)^2}}\leq\frac{\Phi(x)}{\k}.
$$
Using now inequality~\eqref{eq-*} one has
$$
\frac{1}{\Phi'(x)}\geq\frac{1}{\Phi^*\PAR{(1-\ep)\Phi(x)}},
$$
then for all $x\geq M$,
$$
\tau_2\PAR{(1+\ep)^2\frac{\Phi(x)}{\Phi^*\PAR{(1-\ep)\Phi(x)}^2}}\leq\frac{\Phi(x)}{\k}.
$$
Take now $z=(1-\ep) \Phi(x)$,

$$
\tau_2\PAR{\frac{(1+\ep)^2}{1-\ep}\frac{z}{\Phi^*\PAR{z}^2}}\leq\frac{z}{(1-\ep)\k},
$$
to finish take $x=\exp\PAR{\frac{4z}{(1-\ep)\k}}$ to obtain
$$
\tau_2\PAR{{(1+\ep)^2\k}\frac{\log x^2}{\Phi^*\PAR{\frac{(1-\ep)\k}{2}\log x^2}^2}}\leq\frac{1}{2}{\log x^2}.
$$

Recall that $\la=(1-\ep)\k/2$ and let take
$
u_0={1}/\PAR{(1+\ep)^2\k},
$
to obtain the result for
$x\geq C$,
where $C$ is a constant depending on $\Phi$.

If we have $A_\la<C$, one can change the value of $u_0$ to obtain also the results on $[A_\la,C]$.
\end{eproof}

{\noindent {\emph{\textbf{Proof of Theorem~\ref{theo-ls}}}}\\\proofbegin~}
To give the proof of the theorem we need to give an other result like Proposition~\ref{prop-grand2}.
By the same argument as in Proposition~\ref{prop-grand2} one can also prove that there exists
$A,A',B,D>0$ such that for any functions $f\geq0$ satisfying
$$\int f^2 d\mu_\al=1\,\, \text{and} \,\,\ent{\mu_\al}{f^2}\geq 1$$
 we have for some $C'(A_\lambda)$, $C(A_\lambda)$
\begin{equation}
\label{eqfin}
\ent{\mu_\Phi}{f^2}\leq C'(A_\lambda)\var{\mu_\Phi}{f}+C(A_\lambda)\int_{f^2\geq
A_\la}H\PAR{{\frac{f'}{f}}}d\mu_\Phi,
\end{equation}
where $H_\Phi$ is defined on~\eqref{defh} and $A_\la$ on the Proposition~\ref{prop-petit}.
To introduce $A_\la$, we just have to change constants in the inequality.

Then the proof of the theorem is a simple consequence of~\eqref{eqfin} and Proposition~\ref{prop-petit}.
\proofend


\section{Classical properties and applications}
\label{sec-app}
Let us give here properties inherited directly from the methodology known for classical logarithmic Sobolev inequalities.

\begin{eprop}
\label{prop-triv}
\begin{enumerate}[1]
\item This property is known under the name of tensorisation.

Let $\mu_1$ and $\mu_2$ two probability measures on $\dR^{n_1}$
and $\dR^{n_2}$. Suppose that $\mu_1$ (resp. $\mu_2$) satisfies
the a $LSI$ with function $H_\Phi$ and constant $A_1$ (resp. with
constant $A_2$) then the probability $\mu_1\otimes\mu_2$ on
$\dR^{n_1+n_2}$, satisfies a $LSI$ with function $H_\Phi$ and
constant $\max\BRA{A_1,A_2}$.

\item  This property is known under the name of perturbation.

Let $\mu$ a measure on $\dR^n$ a $LSI$ with function $H_\Phi$ and
constant $A$. Let $h$ a bounded function on $\dR^n$ and defined
$\tilde{\mu}$ as
$$
d\tilde{\mu}=\frac{e^h}{{Z}}d\mu,
$$
where $Z=\int e^h d\mu$.

Then the measure $\tilde{\mu}$ satisfies a $LSI$ with function
$H_\Phi$ and the constant $D=Ae^{2\text{osc}(h)}$, where
$\text{osc}(h)=\sup(h)-\inf(h)$.
\item Link between  $LSI$  of function $H_\Phi$ with  Poincar\'e inequality.

Let $\mu$ a measure on  $\dR^n$. If $\mu$  satisfies a $LSI$ with
function $H_\Phi$ and constant $A$, then $\mu$ satisfies a
Poincar\'e inequality with the constant $A$. Let us recall that
$\mu$ satisfies a Poincar\'e  inequality with constant $A$ if
\begin{equation*}
\var{\mu}{f}\leq A\int \ABS{\nabla f}^2d\mu,
\end{equation*}
for all smooth function $f$.
\end{enumerate}

\begin{eproof}
One can find the details of the proof of the properties of tensorisation and  perturbation
 and the implication of the Poincar\'e inequality in chapters 1 and 3  of
\cite{logsob} (Section~1.2.6., Theorem~3.2.1 and Theorem~3.4.3).
\end{eproof}
\end{eprop}

\begin{eprop}
\label{prop-cons}
  Assume that the probability measure $\mu$ on $\dR$ satisfies  a $LSI$ with function $H_\Phi$ and constant $A$. Then there exists
constants $B,C,D\geq0$, independent of $n$ such that:
  if $F$ is a function on $\dR^n$ such that $\forall i$,
$\NRM{\partial_iF}_\infty\leq \zeta$, then we get for $\la\geq0$,
\begin{equation}
\label{eq-cons}
\mu^{\otimes n}(\ABS{F -\mu^{\otimes n}(F)}\geq \la)\leq
\left\{
\begin{array}{ll}
\disp2\exp\PAR{-nB\Phi\PAR{C\frac{\la}{n\zeta}}}&\text{if }
\la > {nD\zeta},\\
\disp2\exp\PAR{-B\frac{\la^2}{n\zeta^2}}& \text{if }0
\leq\la\leq {nD\zeta}.
\end{array}
\right.
\end{equation}
\end{eprop}

\begin{eproof}
Let us first present the proof when $n=1$. Assume, without loss of generality, that  $\int Fd\mu=0$. Due to the homogeneous
property of~\eqref{eq-cons} on can suppose that $\zeta=1$.

Let us recall briefly Herbst's argument (see Chapter~7 \cite{logsob} for more details). Denote $\psi(t)=\int e^{t F}d\mu$,
and remark that $LSI$ of function $H_\Phi$ applied to $f^2=e^{tF}$, using basic properties of $H_\Phi$, yields to
\begin{equation}\label{hehehe}
t\psi'(t)-\psi(t)\log \psi(t)\le A H_\Phi\left(\frac{t}{2}\right)\psi(t)
\end{equation}
which, denoting $K(t)=(1/t)\log \psi(t)$, entails
$$ K'(t)\leq {\frac{A}{t^2}}H_\Phi\PAR{\frac{t}{2}}.$$
Then, integrating, and using $K(0)=\int Fd\mu=0$, we obtain
\begin{equation}\label{laplace}
\psi(t)\le \exp\left( At\int_0^t {1\over s^2}H_\Phi\left(\frac{s}{2}\right)ds\right).
\end{equation}
Then we get using Markov inequality
$$
\mu(\ABS{F -\mu(F)}\geq \la)\leq 2 \exp\PAR{\min_{t\geq0}\BRA{At\int_0^t {\frac{1}{s^2}}H_\Phi\left(\frac{s}{2}\right)ds-\la t}}.
$$
Let note, for $t\geq0$,
$$
G(t)=At\int_0^t {\frac{1}{s^2}}H_\Phi\PAR{\frac{s}{2}}ds-\la t.
$$
An easy study proves that $G$ admits a single minimum on $\dR^+$ (except maybe
 if $\lambda=0$).
Then due to the definition of $H_\Phi$  we get that
$$
\min_{t\geq0}\BRA{G(t)}=-\frac{\la^2}{A},\quad\text{if }\la\leq
AD.
$$
Assume now that $\la\geq AD$ then  we obtain after derivation
\begin{equation}
\label{eq-gr}
\min_{t\geq0}\BRA{G(t)}=-A\Phi^*\PAR{t_0\frac{B}{2}},\quad \text{with }\la t_0=At_0\int_0^{t_0} {\frac{1}{s^2}}
H_\Phi\PAR{\frac{s}{2}}ds+AH_\Phi\PAR{\frac{t_0}{2}}.
\end{equation}
We first prove that there exists $C\geq0$ such that  for all $t_0$
large enough
\begin{equation}
\label{eq-fi}
t_0\int_0^{t_0} {\frac{1}{s^2}}H_\Phi\PAR{\frac{s}{2}}ds\leq CH_\Phi\PAR{\frac{t_0}{2}}.
\end{equation}
For
$\k\geq0$ large enough and $t_0\geq \k$ we get using then inequality~\eqref{eq-*} we get
$$
t_0\int_\k^{t_0} {\frac{1}{s^2}}H_\Phi\PAR{\frac{s}{2}}ds
\leq Ct_0\int_\k^{t_0} {\frac{1}{s^2}}\Phi\PAR{\Phi'^{-1}\PAR{\frac{s}{2}}}ds,
$$
with $C\geq0$.  Then by a change of variables and integration by
parts, for large enough $t_0$,
$$
\begin{array}{rl}
\disp t_0\int_\k^{t_0} {\frac{1}{s^2}}\Phi\PAR{\Phi'^{-1}\PAR{\frac{s}{2}}}ds
\disp & =\disp\frac{t_0}{2}\int^{ \Phi'^{-1}\PAR{\frac{t_0}{2}}}_{ \Phi'^{-1}\PAR{\frac{\k}{2}}}
\frac{\Phi(u)}{\Phi'(u)^2}\Phi''(u)du\\
\disp &\disp \leq \frac{t_0}{2}\frac{\Phi(\Phi^{-1}(\k/2))}{\Phi'(\Phi^{-1}(\k/2))}
+\frac{t_0}{2}{\Phi'}^{-1}(t_0/2)\\
\disp &\disp\leq Ct_0{\Phi'}^{-1}(t_0/2),
\end{array}
$$
for some other $C\geq0$. Then we get, using
inequality~\eqref{eq-lem3}, for $t_0$ large enough,
$$
t_0\int_0^{t_0} {\frac{1}{s^2}}H_\Phi\PAR{\frac{s}{2}}ds\leq
Ct_0{\Phi'}^{-1}(t_0/2)\leq C'\Phi^*(t_0/2).
$$
for some constant $C'\geq0$ and for $t_0$ large enough and inequality~\eqref{eq-fi} is proved. By~\eqref{eq-fi}
and~\eqref{eq-gr} one get for $t_0$ large enough,
$$
\la t_0\leq A'\Phi^*\PAR{\frac{t_0}{2}},
$$
for some constant $A'\geq0$. But, using  inequality~\eqref{eq-lem3} we get then
$$
\Phi'(A\la)\leq Ct_0,
$$
$$
\min_{t\geq0}\BRA{G(t)}\leq -A\Phi^*\PAR{B\Phi'(C\la)}\leq -A\Phi^*\PAR{\Phi'(C'\la)},
$$
if $\la$ is large enough and for some
other constants $A,B,C,C'\geq0$. Using
 inequality~\eqref{eq-*},
we obtain the result in dimension 1.

\medskip

For the $n$-dimensional extension, use the tensorisation property of $LSI$ of function $H_\Phi$ and
$$
\sum_{i=1}^n H_\Phi\PAR{\frac{t}{2}\partial_i F}\leq nH_\Phi\PAR{\frac{t}{2}}.
$$
Then we can use the case of dimension 1 with the constant $A$ replaced by $An$.
\end{eproof}

\begin{erem}
Let us present a simple application of the preceding
 proposition
to deviation inequality of the empirical mean of a function.
Consider the real valued function $f$, with $|f'|\le 1$. Let apply
Proposition~\ref{prop-cons} with the two
 functions
$$F(x_1,...,x_n)={1\over n}\sum_{k=1}^nf(x_k)\quad \text{and}\quad
F(x_1,...,x_n)={1\over\sqrt{n}}\sum_{k=1}^n f(x_i).$$ We obtain
then
\begin{equation*}
\dP\left({1\over n}\ABS{\sum_{k=1}^nf(X_k)-\mu(f)}>\la\right)\leq
\left\{
\begin{array}{ll}
\disp 2\exp\PAR{-nA \Phi(B\la)}&\text{if }
\la\geq D,\\
\disp 2\exp\PAR{-nA{\la^2}}& \text{if }0\leq\la\leq D,
\end{array}
\right.
\end{equation*}
\begin{equation*}
\dP\left({1\over
\sqrt{n}}\ABS{\sum_{k=1}^nf(X_k)-\mu(f)}>\la\right)\leq \left\{
\begin{array}{ll}
\disp 2\exp\PAR{-nA \Phi\PAR{B\frac{\la}{\sqrt{n}}}}&\text{if }
\la\geq D\sqrt{n},\\
\disp 2\exp\PAR{-A{\la^2}}& \text{if }0\leq\la\leq D\sqrt{n}.
\end{array}
\right.
\end{equation*}
\end{erem}


\newcommand{\etalchar}[1]{$^{#1}$}

Ivan Gentil, Arnaud Guillin\\
CEREMADE \\
(UMR 7534, Universit\'e Paris-Dauphine et CNRS)\\
Place du Mar\'echal De Lattre de Tassigny, 75775 Paris C\'edex~16, France\\
\{gentil,guillin\}@ceremade.dauphine.fr\\
http://www.ceremade.dauphine.fr/
\raisebox{-4pt}{$\!\!\widetilde{\phantom{x}}$}\{gentil,guillin\}\\

Laurent Miclo\\
Laboratoire d'Analyse, Topologie, et Probabilit\'es \\
(UMR  6632, Université de Provence et CNRS) \\
39, rue F. Joliot Curie, 13453 Marseille C\'edex 13, France\\
miclo@cmi.univ-mrs.fr

\end{document}